\documentclass[11pt]{article}
\usepackage{amsmath}
\usepackage{amssymb}
\usepackage{amsfonts}
\usepackage{amscd}
\usepackage{enumerate}
\usepackage{fancyheadings}
\usepackage{color}
\usepackage{graphicx}
\usepackage{amsthm}
\newtheorem{Lem}{Lemma}

\newtheorem{Th}{Theorem}
\newtheorem{Co}{Corollary}
\newtheorem{Def}{Definition}

\DeclareMathOperator{\J}{J}

\DeclareMathOperator{\supp}{supp}

\DeclareMathOperator{\graph}{graph}
\DeclareMathOperator{\esupp}{ess-supp}
\DeclareMathOperator{\diag}{diag}

\DeclareMathOperator{\df}{def}

\begin{document}

\begin{center}
\bf{\Large Structure of the Short Range Amplitude for General Scattering Relations}
\end{center}

\noindent
{\bf Ivana Alexandrova}

\noindent
Department of Mathematics, University of Toronto, Toronto, Ontario, Canada
M5S 3G3, Tel.: 1-416-946-0318, Fax: 1-416-978-4107, email: 
alexandr@math.toronto.edu

\noindent
November 25, 2004

\begin{abstract}
We consider scattering by short range perturbations of the 
semi-classical Laplacian. 
We prove that when a polynomial bound on the resolvent holds the 
scattering amplitude is a semi-classical Fourier integral operator associated 
to the scattering relation near a non-trapped ray.
Compared to previous work, we allow the scattering relation to have more 
general structure.
\end{abstract}

\noindent
{\bf Keywords and phrases:} Short range perturbations, scattering 
amplitude, scattering relation, 
semi-classical Fourier integral operators.

\section{Introduction and Statement of Results}
We study the structure of the scattering amplitude associated with the 
semi-classical Schr\"{o}dinger operator with a short range potential on $\mathbb{R}^{n}.$
We prove that, when restricted away from the diagonal on 
$\mathbb{S}^{n}\times\mathbb{S}^{n},$ the natural scattering amplitude quantizes the scattering 
relation in the sense of semi-classical Fourier integral operators.
The scattering relation at energy $\lambda>0$ here is given roughly by the Hamiltonian flow of the symbol $p$ of the operator between two hypersurfaces ``at infinity'' inside the energy surface $\{p=\lambda\}.$

\subsection{A Survey of Earlier Results}

The structure of the scattering matrix has been of significant interest 
to researchers in mathematical physics.
Earlier results have focused primarily on establishing asymptotic expansions of the 
scattering amplitude.
In this section we describe briefly only those asymptotic expansions most relevant to our 
work and refer to \cite{AI} for a more comprehensive survey.

We begin by introducing some notation.
Let $P(h)=-\frac{1}{2}h^{2}\Delta+V,$ $0<h<<1,$ where
\begin{equation}\label{potential}
\left|\partial^{\alpha}V(x)\right|\leq 
C_{\alpha}\langle x\rangle^{-\rho-|\alpha|},\; 
x\in\mathbb{R}^{n}, 
\rho>1,
\end{equation}
where $\langle x\rangle=(1+\|x\|^{2})^\frac{1}{2}.$
Let $\lambda>0$ and for $\omega\in\mathbb{S}^{n-1}$ and 
$z\in\omega^{\perp}$ we denote by
\[\gamma_{\infty}\left(\cdot; z, \sqrt{2\lambda}\omega\right)=\left\{q_{\infty}\left(\cdot; z, \sqrt{2\lambda}\omega\right), 
p_{\infty}\left(\cdot; z, \sqrt{2\lambda}\omega\right)\right\}\] 
the unique phase 
trajectory, i.~e. the integral curve of the Hamiltonian vector field of 
$p(x, \xi)=\frac{1}{2}\|\xi\|^{2}+V(x),$ such that 
\begin{gather*}
\lim_{t\to -\infty}\left\|q_{\infty}\left(t; z, 
\sqrt{2\lambda}\omega\right)-\sqrt{2\lambda}\omega 
t-z\right\|=0,\\
\lim_{t\to -\infty}\left\|p_{\infty}\left(t; z, \sqrt{2\lambda}\omega\right)-\sqrt{2\lambda}\omega\right\|=0
\end{gather*}
in the $C^{\infty}$ topology for the impact parameter $z$ and $\omega.$
If $\lim_{t\to\infty}\left\|q_{\infty}\left(t; z, \sqrt{2\lambda}\omega\right)\right\|=\infty,$ then, setting $\mathbb{S}^{n-1}_{2\lambda}=\{x\in\mathbb{R}^{n}: \|x\|=2\lambda\},$ we have that there exist $U\subset T^{\star}\mathbb{S}^{n-1}_{2\lambda}$ open, $\left(\sqrt{2\lambda}\omega, z\right)\in U,$ where $\mathbb{S}^{n-1}_{2\lambda}=\{x\in\mathbb{R}^{n}: \|x\|=2\lambda\},$ $\xi_{\infty}\in C^{\infty}\left(T^{*}\mathbb{S}^{n-1}_{2\lambda}\cap U; \mathbb{S}^{n-1}\right),$
and $x_{\infty}\in C^{\infty}\left(T^{*}\mathbb{S}^{n-1}_{2\lambda}\cap U; \mathbb{R}^{n}\right)$ 
such that 
\begin{gather*}
\lim_{t\to\infty}\left\|q_{\infty}\left(t; z, 
\sqrt{2\lambda}\omega\right)-\sqrt{2\lambda}\xi_{\infty}\left(z, 
\sqrt{2\lambda}\omega\right)t-x_{\infty}\left(z, 
\sqrt{2\lambda}\omega\right)\right\|_{C^{\infty}(U)}=0,\\
\lim_{t\to\infty}\left\|q_{\infty}\left(t; 
z, \sqrt{2\lambda}\omega\right)-\sqrt{2\lambda}\xi_{\infty}\left(z, 
\sqrt{2\lambda}\omega\right)\right\|_{C^{\infty}(U)}=0.
\end{gather*}
The trajectory $\gamma_{\infty}\left(\cdot; z, \sqrt{2\lambda}\omega\right)$ is
then said to have initial 
direction $\omega$ and final direction $\theta=\xi_{\infty}\left(z, 
\sqrt{2\lambda}\omega\right).$
We also make the following
\begin{Def}
The outgoing direction $\theta\in\mathbb{S}^{n-1}$ is called {\it non-degenerate}, or 
{\it regular}, for 
the incoming direction $\omega\in\mathbb{S}^{n-1}$ if $\theta\ne\omega$ and for all 
$z'\in\omega^{\perp}$ with 
$\xi_{\infty}\left(z', \sqrt{2\lambda}\omega\right)=\theta,$ the map $\omega^{\perp}\ni z\mapsto \xi_{\infty}\left(z, \sqrt{2\lambda}\omega\right)\in\mathbb{S}^{n-1}$ is non-degenerate at $z'.$
\end{Def}

Several authors, working under the assumption that a certain final direction $\theta$ is non-degenerate for a given initial direction $\omega,$ have proved asymptotic expansions of the scattering amplitude $A$ of the form
\begin{equation}\label{vexpansion}
K_{A(\lambda, h)}(\omega, \theta)=\sum_{j=1}^{l}\hat{\sigma}\left(z_{j}, \omega; 
\lambda\right)^{-1/2}\exp\left(ih^{-1}S_{j}-i\mu_{j}\pi/2\right)+\mathcal{O}\left(h\right),
\end{equation}
where $\left(z_{j}\right)_{j=1}^{l}\equiv\left(\xi_{\infty}^{-1}\left(\cdot, 
\sqrt{2\lambda}\omega\right)\right)
(\theta_0),$  $\hat{\sigma}\left(z_{j}, \omega; 
\lambda\right)=\det\left(\J\xi_{\infty}\left(\cdot, \sqrt{2\lambda}\omega\right)\right)\left(z_j\right),$ with $\J$ denoting the Jacobian matrix, 
\begin{equation}\label{modaction}
S_j=\int_{-\infty}^{\infty}\left(\frac{1}{2}\left|p_{\infty}\left(t; z, \sqrt{2\lambda}\omega\right)\right|^{2}-V\left(q_{\infty}\left(t; z, 
\sqrt{2\lambda}\omega\right)\right)-\lambda\right)dt-\left\langle x_{\infty}\left(z, \sqrt{2\lambda}\omega\right),
\sqrt{2\lambda}\theta\right\rangle
\end{equation}
is a modified action along the $j-$th $(\omega, \theta)$ trajectory, 
and $\mu_{j}$ is the path index of that trajectory.
Vainberg \cite{V} has studied smooth compactly supported potentials $V$ at energies 
$\lambda>\sup V$ and has proved such an asymptotic expansion with the error term 
estimated uniformly over a sufficiently small neighborhood containing the final direction 
while the initial direction is held a constant.
Guillemin \cite{G} has established a similar asymptotic expansion in the setting of smooth 
compactly-supported metric perturbations of the Laplacian for fixed incoming and outgoing 
directions.
Working with non-trapping potential perturbations of the Laplacian satisfying 
\eqref{potential} with $\rho>\max\left(1, \frac{n-1}{2}\right),$ Yajima 
\cite{Y} has proved such an asymptotic expansion in the $L^{2}$ sense.
For non-trapping short-range ($\rho>1$) potential perturbations of the Laplacian, Robert and 
Tamura \cite{RT} have established a pointwise asymptotic expansion of this form for 
constant initial and final directions.
This result has been extended to the case of 
trapping energies by 
Michel \cite{M} under an additional assumption on the distribution of the resonances of $P(h).$

In \cite{AI} we have proved, without making the non-degeneracy assumption, that the scattering 
amplitude for smooth compactly supported potential and metric perturbations of the Euclidean 
Laplacian at both trapping and non-trapping energies is a semi-classical Fourier integral 
operator associated to the scattering relation.
We have further showed how the expansion \eqref{vexpansion} follows from the general 
theory of semi-classical Fourier integral operators developed in \cite{Afio}, 
once the non-degeneracy assumption on the initial and final directions is made.
Here we extend these results to the case of short-range perturbations of the Laplacian 
when the scattering amplitude is restricted away from the diagonal in 
$\mathbb{S}^{n-1}\times\mathbb{S}^{n-1}.$ 

\subsection{Statement of Main Theorem}
We consider the semi-classical Schr\"{o}dinger operator 
$P(h)=-\frac{1}{2}h^{2}\Delta+V,$ on $\mathbb{R}^{n}$ for $n\geq 2,$ 
$0<h\leq1,$  with the potential $V\in C^{\infty}(\mathbb{R}^{n}; \mathbb{R})$ satisfying \eqref{potential}.
Let $P_{0}(h)=-\frac{1}{2}h^{2}\Delta.$
Then $P(h)$ and $P_0(h)$ admit unique self-adjoint realizations 
on $L^{2}(\mathbb{R}^{n})$ with common domain $H^2(\mathbb{R}^{n}).$ 
It is well-known that the wave operators
\begin{equation*}
W_{\pm}=\text{s-}\lim_{t\to\pm\infty}U(t)U_0(-t),
\end{equation*} 
where 
\[U(t)=e^{-\frac{i}{h}tP(h)}, \; U_0(t)=e^{-\frac{i}{h}tP_{0}(h)},\; t\in\mathbb{R}.\]
We can therefore define the scattering operator
\begin{equation*}
S=W_{+}^{*}W_{-}=\mathcal{F}^{-1}\int_{\lambda>0}\bigoplus S(\lambda, h) 
d\lambda\,\mathcal{F},
\end{equation*}
where $\mathcal{F}$ denotes the unitary Fourier transform on $L^{2}(\mathbb{R}^{n}).$
The operator $S(\lambda, h)$ is called the scattering matrix at energy 
$\lambda>0$ and is a unitary operator on 
$L^{2}(\mathbb{S}^{n-1}).$ 
The scattering amplitude $A(\lambda, h)$ is defined by 
$A(\lambda, h)=c(n, \lambda, h)T(\lambda, h),$ where $T(\lambda, 
h)=-i(2\pi)^{-1}\left(I-S(\lambda, h)\right)$ and 
\[c(n, \lambda, h)=-2\pi(2\lambda)^{-\frac{n-1}{4}}(2\pi 
h)^{\frac{n-1}{2}}e^{-i\frac{(n-3)\pi}{4}}.\]

To state our Main Theorem, we let $R(\lambda+i0, h)=\lim_{\epsilon\downarrow
0}\left(P(h)-\lambda-i\epsilon\right)^{-1},$ where the limit is taken in
the space $\mathcal{B}(L^{2}_{\alpha}(\mathbb{R}^{n}),
L^{2}_{-\alpha}(\mathbb{R}^{n})),$
$\alpha>\frac{1}{2},$ with $L^{2}_{\alpha}(\mathbb{R}^{n})=\{f:
\langle\cdot\rangle^{\alpha}f\in L^{2}(\mathbb{R}^{n})\}.$
We further refer the reader to Section \ref{sgeom} for the definitions non-trapped trajectories and the scattering relation $SR_{\bar{U}}(\lambda).$ 
The class of semi-classical Fourier integral operators $\mathcal{I}^{r}_{h}$ is defined in Appendix \ref{scanal},
where we also review the notion of pseudodifferential operators of principal type.   

We are now ready to prove our 

\medskip
\noindent
{\bf Main Theorem.}
{\it Let $\lambda>0$ be such that the operator $P(h)-\lambda$ is of principal type.
Let also 
\begin{equation}\label{resest}
\left\|R(\lambda+i 0, h)\right\|_{\mathcal{B}\left(L^{2}_{\alpha}(\mathbb{R}^{n}), 
L^{2}_{-\alpha}(\mathbb{R}^{n})\right)}=\mathcal{O}(h^{s}), \;s\in\mathbb{R}, \;\alpha>\frac{1}{2}.
\end{equation} 
Let $(\omega, z)\in T^{*}\mathbb{S}^{n-1}$ be such that 
$\gamma_{\infty}\left(\cdot; z, \sqrt{2\lambda}\omega\right)$ is a 
non-trapped trajectory.

Then there exist open sets $U\subset T^{*}\mathbb{S}^{n-1}$ with 
$(\omega, z)\in U$ such that 
{\em 
\[A(\lambda, h)\in \mathcal{I}_{h}^{\frac{n}{2}+2}\left({\mathbb
S}^{n-1} \times{\mathbb S}^{n-1}\backslash\diag(\mathbb{S}^{n-1}\times\mathbb{S}^{n-1}), 
SR_{\bar{U}}(\lambda)\right).\]}}

We remark that \cite[Theorem 1]{A} gives a characterization of semi-classical Fourier integral distributions as oscillatory integrals.
Applied to the scattering amplitude here this characterization says approximately that for every non-degenerate phase function $\phi$ which locally parameterizes $SR_{\bar{U}}(\lambda)$ we can find a symbol $a$ admitting an asymptotic expansion in $h$ such that $CK_{A(\lambda, h)},$ where $C$ is a microlocal cut-off to $SR_{\bar{U}}(\lambda)$ (see Appendix \ref{scanal}), can be represented as an oscillatory integral with phase $\phi$ and symbol $a.$ 
From the discussion in \cite[Section 4.1]{A} we further know that such a non-degenerate 
phase function always exists, and therefore we can always express $CK_{A(\lambda, h)}$ as an oscillatory integral admitting an asymptotic expansion in $h.$ 
In the special case when the non-degeneracy assumption holds, we recover the phases \eqref{maction} in \eqref{vexpansion} -- see 
Theorem \ref{tmicrol} below. 
We expect that a finer analysis based on our method would give a precise description of the 
amplitudes as well. 
What is different here is the fact that we can handle the cases in which the non-degeneracy 
assumption fails.

We now introduce some of the notation we shall use below.
For a 
sequentially continuous operator
$T:C^{\infty}_{c}(\mathbb{R}^{m})\to\mathcal{D}'(\mathbb{R}^{n})$ we shall
denote by $K_{T}$ its Schwartz kernel.
On any smooth manifold $M$ we denote by $\sigma$ the canonical symplectic
form on $T^{*}M$ and everywhere below we work with the canonical
symplectic structure on $T^{*}M.$
We shall denote by $H_{p}$ the Hamiltonian vector field of $p.$
The integral curve of $H_{p}$ with initial conditions $(x_0, \xi_0)\in 
T^{*}\mathbb{R}^{n}$ will
be denoted by
$\gamma(\cdot; x_0, \xi_0)=(x\left(\cdot; x_0, \xi_0\right), \xi(\cdot;
x_0, \xi_0)).$
If $C\subset T^{*}M_1\times T^{*}M_2,$ where $M_j,$ $j=1, 2,$ are smooth
manifolds, we will use the notation $C'=\{(x, \xi; y,
-\eta): (x, \xi; y, \eta)\in C\}.$
We shall also use $\|\cdot\|_{\pm\gamma, \mp\gamma}$ to denote the norm of a linear operator 
between the spaces $L^{2}_{\pm\gamma}(\mathbb{R}^{n})$ and 
$L^{2}_{\mp\gamma}(\mathbb{R}^{n}).$
Lastly, we set $B(0, r)=\{x\in\mathbb{R}^{n}: \|x\|<r\}$ and $B(0, r, 
r+1)=\{x\in\mathbb{R}^{n}: r<\|x\|<r+1\}.$

This paper is organized as follows. 
We review the definition of semi-classical Fourier integral distributions and 
operators in the Appendix, where we also recall the relevant part of 
semi-classical analysis. 
Isozaki-Kitada's representation of the short-range scattering amplitude which 
we will use in this article is presented in Section \ref{sreprsa}.
Two preliminary lemmas giving additional information on the structure of 
the semi-classical amplitude, are given in Section \ref{2ls}.
The scattering relation is defined in Section \ref{sgeom}, where we also prove that it can 
be parameterize by the modified actions when the non-degeneracy assumption is made.
The proof of the Main Theorem is presented in Section \ref{pmain} and its applications to 
non-trapping and trapping perturbations are discussed in Section \ref{snt} and Section 
\ref{str}, respectively.
Finally, the theorem giving the microlocal representation of the scattering amplitude 
as an oscillatory integral under 
the 
non-degeneracy assumption is proved in Section \ref{smicrol}.

\section{Preliminaries}
In this section we introduce some of the preliminary results we shall use throughout this article.

\subsection{Representation of the Scattering Amplitude}\label{sreprsa}
Here we present the representation of the short range scattering amplitude 
developed by Isozaki and Kitada \cite{IK}.
This is the representation we shall use in the proof of our Main 
Theorem.

\begin{Def}
Let $\Omega\subset T^{*}\mathbb{R}^{n}$ be an open subset.
We denote by $A_{m}(\Omega)$ the class of symbols $a$ such that $(x, \xi)\mapsto a(x, \xi, 
h)$ belongs to $C^{\infty}(\Omega)$ and
\begin{equation*}
\left|\partial_{x}^{\alpha}\partial_{\xi}^{\beta}a(x, \xi)\right|\leq C_{\alpha\beta}\langle 
x \rangle^{m-|\alpha|}\langle\xi\rangle^{-L}, \text{ for all } (x, \xi)\in\Omega, (\alpha, 
\beta)\in\mathbb{N}^{n}\times\mathbb{N}^{n}, L>0.
\end{equation*}
We denote $A_{m}(T^{*}\mathbb{R}^{n})$ by $A_{m}.$ 
\end{Def}

We also use the notation 
\[\Gamma_{\pm}(R, d, \sigma)=\left\{(x, \xi)\in\mathbb{R}^{n}\times\mathbb{R}^{n}:|x|>R, 
\frac{1}{d}<|\xi|<d, \pm cos(x, \xi)>\pm\sigma\right\}\] with $R>1,$ 
$d>1,$ $\sigma\in(-1, 
1),$ and $cos(x, \xi)=\frac{\langle x, \xi\rangle}{|x||\xi|},$ for the 
outgoing and incoming 
subsets of phase space, respectively.

For $\alpha>\frac{1}{2},$ we denote the bounded operator $F_{0}(\lambda, h):
L^{2}_{\alpha}(\mathbb{R}^{n})\rightarrow L^{2}(\mathbb{S}^{n-1})$ by
\begin{equation*}
\left(F_{0}(\lambda, h)f\right)(\omega)=(2\pi 
h)^{-\frac{n}{2}}(2\lambda)^{\frac{n-2}{4}}\int_{\mathbb{R}^{n}}
e^{-\frac{i}{h}\sqrt{2\lambda}\langle \omega, x\rangle}f(x)dx, \lambda>0.
\end{equation*}

Let $R_{0}>>0,$ $1<d_4<d_3<d_2<d_1<d_0,$ and 
$0<\sigma_4<\sigma_3<\sigma_2<\sigma_1<\sigma_0<1.$
Using the WKB method, Isozaki and Kitada \cite{IK} have constructed 
parametrices for the wave operators with phase functions 
$\Phi_{\pm}$ and symbols $a_{\pm}$ and $b_{\pm}$ such that:
\begin{enumerate}
\item $\Phi_{\pm}\in C^{\infty}(T^{*}\mathbb{R}^{n})$ solve the eikonal equation 
\begin{equation}\label{eikonal}
\frac{1}{2}\left|\nabla_{x}\Phi_{\pm}(x, \xi)\right|^{2}+V(x)=\frac{1}{2}|\xi|^{2}
\end{equation}
in $(x, \xi)\in\Gamma_{\pm}(R_0, d_0, \pm\sigma_0),$ respectively.
\item $\Phi_{\pm}(\cdot, \cdot\cdot)-\langle\cdot, \cdot\cdot\rangle\in 
A_{0}\left(\Gamma_{\pm}(R_0, d_0, \pm\sigma_{0})\right).$ 
\item For all $(x, \xi)\in T^{*}\mathbb{R}^{n}$
\begin{equation}\label{phinondeg}
\left|\frac{\partial^{2}\Phi_{\pm}}{\partial_{x_{j}}\partial_{\xi_{k}}}(x, 
\xi)-\delta_{jk}\right|<\epsilon(R_0), 
\end{equation}
where $\delta_{jk}$ is the Kronecker delta and $\epsilon(R_0)\to R_{0}$ as $R_{0}\to\infty.$
\item $a_{\pm}\sim\sum_{j=0}^{\infty}h^j a_{\pm j},$ where $a_{\pm j}\in 
A_{-j}(\Gamma_{\pm}(3R_0, d_1, \mp\sigma_1)),$ $\supp a_{\pm j}\subset\Gamma_{\pm}(3R_0, 
d_1, \mp\sigma_1),$ $a_{\pm j}$ solve
\begin{equation}\label{aeq1}
\langle \nabla_{x}\Phi_{\pm}, \nabla_{x}a_{\pm 0}\rangle 
+\frac{1}{2}\left(\Delta_{x}\Phi_{\pm}\right)a_{\pm 0}=0
\end{equation}
\begin{equation}\label{aeq2}
\langle \nabla_{x}\Phi_{\pm}, \nabla_{x}a_{\pm 
j}\rangle+\frac12\left(\Delta_{x}\Phi_{\pm}\right)a_{\pm j}=\frac{i}{2} \Delta_{x}a_{\pm 
j-1}, j\geq 1, 
\end{equation}
with the conditions at infinity
\begin{equation}\label{condinfty}
a_{\pm 0}\to 1, a_{\pm j}\to 0, j\geq 1, \text{ as } |x|\to \infty.
\end{equation}
in $\Gamma_{\pm}(2R_0, d_2, \mp\sigma_2),$
and solve \eqref{aeq1} and \eqref{aeq2} in $\Gamma_{\pm}(4R_0, d_1, \mp\sigma_2).$
\item $b_{\pm}\sim\sum_{j=0}^{\infty}h^j b_{\pm j},$ where $b_{\pm j}\in
A_{-j}(\Gamma_{\pm}(5R_0, d_3, \pm\sigma_{4}),$ $\supp b_{\pm j}\subset\Gamma_{\pm}(5R_0, 
d_3, \pm\sigma_{4}),$ $b_{\pm j}$ solve \eqref{aeq1} and \eqref{aeq2} with the conditions 
at infinity \eqref{condinfty} in $\Gamma_{\pm}(6R_0, d_4, \pm\sigma_3),$ and solve 
\eqref{aeq1} and \eqref{aeq2} in $\Gamma_{\pm}(6R_0, d_3, \pm\sigma_{3}).$
\end{enumerate}

For a symbol $c$ and a phase function $\phi$, we denote by $I_{h}(c, \phi)$ the oscillatory integral 
\begin{equation*}
I_{h}(c, \phi)=\frac{1}{(2\pi h)^{n}}\int_{\mathbb{R}^{n}} e^{\frac{i}{h}(\phi(x, \xi)-\langle y, 
\xi\rangle)}c(x, \xi) d\xi
\end{equation*}
and let 
\begin{equation*}
\begin{aligned}
K_{\pm a}(h) & =P(h)I_{h}(a_{\pm}, \Phi_{\pm})-I_{h}(a_{\pm}, \Phi_{\pm})P_{0}(h)\\
K_{\pm b}(h) & =P(h)I_{h}(b_{\pm}, \Phi_{\pm})-I_{h}(b_{\pm}, \Phi_{\pm})P_{0}(h).
\end{aligned}
\end{equation*} 

The operator $T(\lambda, h)$ for $\lambda\in\left(\frac{1}{2d_{4}^{2}},
\frac{d_{4}^{2}}{2}\right)$ is then given by (see \cite[Theorem 3.3]{IK})
\begin{equation*}
T(\lambda, h)=T_{+1}(\lambda, h)+T_{-1}(\lambda, h)-T_{2}(\lambda, h),
\end{equation*}
where
\begin{equation*}
T_{\pm 1}(\lambda, h)=F_{0}(\lambda, h)I_{h}(a_{\pm}, \Phi_{\pm})^{*}K_{\pm b}(h)
F_{0}^{*}(\lambda, h)
\end{equation*}
and
\begin{equation*}
T_{2}(\lambda, h)=F_{0}(\lambda, 
h)K_{+a}^{*}(h)R(\lambda+i0, 
h)\left(K_{+b}(h)+K_{-b}(h)\right)F_{0}^{*}(\lambda, h),
\end{equation*}

\subsection{Two Preparatory Lemmas}\label{2ls}
The following two lemmas will be useful in studying the structure of the scattering amplitude.  
\begin{Lem}\label{tempop}
Let $W=\mathcal{O}_{\mathcal{B}(L^{2}(\mathbb{R}^{n}))}(h^{s}),$ $h\to 
0,$ or $W=\mathcal{O}_{\mathcal{B}(L^{2}_{\alpha}(\mathbb{R}^{n}), 
L^{2}_{-\alpha}(\mathbb{R}^{n}))}(h^{s}),$ $h\to 0,$ for some $s\in\mathbb{R}.$

Then $K_W\in\mathcal{D}_{h}'(\mathbb{R}^{2n}).$  
\end{Lem}

\begin{proof}
By Schwartz Kernel Theorem, for some $h_0>0$ and every $h\in(0, h_0],$ there 
exists $w_h\in\mathcal{D}'(\mathbb{R}^{2n})$ such that $\langle T \varphi, 
\psi\rangle=\langle w_h, \varphi\otimes\psi\rangle,$ $\varphi, \psi\in 
C_{c}^{\infty}(\mathbb{R}^{n}).$
Let $\chi\in C_{c}^{\infty}(\mathbb{R}^{2n})$ and let $c_1>c_2>0$ be such that $\supp\chi\subset K_1(c_2)\times K_2(c_2),$ where $K_j(d)=\{x\in\mathbb{R}^{n}: |x_l|<d, l=1, 
\dots, n\},$ $j=1, 2,$ $d>0.$
Let also $\rho_j\in 
C_{c}^{\infty}(K_j(c_1)),$ $j=1, 2,$ be such that $\rho_1\times\rho_2=1$ on $K_1(c_2)\times K_2(c_2).$
Then, by the proof of Schwartz Kernel Theorem \cite[Theorem 6.1.1]{F}, we have that
\begin{equation*}
\left\langle w_h, \chi e^{-\frac{i}{h}\left(\langle \cdot, 
\xi\rangle+\langle\cdot\cdot, 
\eta\rangle
\right)}\right\rangle=\Sigma_{\mathbb{Z}^{n}\times\mathbb{Z}^{n}}\hat{\chi}_{m, k}\left\langle T\rho_1 
E_{h}(\langle m, \cdot \rangle), \rho_2 E_{h}(\langle k, \cdot\cdot \rangle)\right\rangle,
\end{equation*}
where $E(t)=e^{\frac{2\pi i t}{b}},$ $t\in\mathbb{R},$ and $\hat{\chi}_{m, 
k}=\frac{1}{b^{2n}}\int_{K_{1}\times K_{2}}\chi(x, y)e^{-\frac{i}{h}\left(\langle x, \xi\rangle+\langle y, \eta\rangle\right)}E(-m\cdot x- k\cdot y) 
dx dy.$
Integration by parts now gives
\begin{equation}\label{est1}
(1+|m|)^{M}(1+|k|)^{M}\hat{\chi}_{m, k}\leq C_{1}h^{-2M}\langle (\xi, 
\eta)\rangle^{M}\sum_{|\alpha|\leq M, |\beta|\leq M}\left\|\partial_{x}^{\alpha}\partial_{y}^{\beta}\chi\right\|_{L^{\infty}(\mathbb{R}^{2n})}, m, k\in \mathbb{Z}^{n}, 
M\in\mathbb{N}_{0}. 
\end{equation}
We also have 
\begin{equation}\label{est2}
\left|\left\langle T \rho_1 E_{h}(\langle m, \cdot \rangle), \rho_2 E_{h}(\langle k, 
\cdot\cdot 
\rangle)\right\rangle\right|\leq C_{2} h^{s}.
\end{equation}
From estimates \eqref{est1} and \eqref{est2} we obtain 
\begin{equation}\label{kernelestim}
\left|\Sigma_{\mathbb{Z}^{n}\times\mathbb{Z}^{n}}\hat{\chi}_{m, k}\left\langle T\rho_1
E_{h}(\langle m, \cdot \rangle), \rho_2 E_{h}(\langle k, \cdot\cdot \rangle)\right\rangle 
\right|\leq C_{3}\sum_{|\alpha|\leq M, |\beta|\leq M}\left 
\|\partial_{x}^{\alpha}\partial_{x}^{\beta}\chi\right\|_{L^{\infty}(\mathbb{R}^{2n})} h^{s-2M}\langle (\xi,
\eta)\rangle^{M},
\end{equation}
with
\begin{equation*}
C_{3}=C_1 C_2 \sum_{\mathbb{Z}^{n}\times\mathbb{Z}^{n}}(1+|m|)^{-M}(1+|k|)^{-M}<\infty,
\end{equation*}
if $M$ is taken large enough.
Therefore $K_T\in\mathcal{D}'_{h}(\mathbb{R}^{2n}).$
\end{proof}

\begin{Lem}\label{rsym}
Let $\nu:\mathbb{R}^{2n}\to\mathbb{R}^{2n}$ be given by $\nu(x, y)=(y, 
x).$

Then $\nu^{*}K_{R(\lambda+i0, h)}=K_{R(\lambda+i0, h)}$ for every 
$\lambda>0.$
\end{Lem}

\begin{proof}
For $u, v\in L^{2}(\mathbb{R}^{n})$ let $\left\langle u, v\right\rangle=\int u v.$
Let $u$ and $v$ further satisfy $u, v\in C_{c}^{\infty}(\mathbb{R}^{n})$ and
let $z\in\mathbb{C}$ be such that $\Im z>0.$
We then have
\begin{equation}\label{symR}
\begin{aligned}
\left\langle R\left(z, h\right)u, v\right\rangle
& =\left\langle R\left(z, h\right)u, \left(P\left(h\right)-z\right)R\left(z,
h\right)v\right\rangle\\
&=\left\langle \left(P(h)-z\right)R\left(z, 
h\right)u, R\left(z, h\right)v\right\rangle\\
&=\left\langle u, R\left(z, h\right)v\right\rangle.
\end{aligned}
\end{equation}

Let, now, $\lambda\in\mathbb{R}\backslash\left\{0\right\}$ and let
$\left(z_{k}\right)_{k\in\mathbb{N}}\subset\mathbb{C}$ satisfy $\Im z_{k}\downarrow 0,$ 
$k\to\infty,$ and $\Re z_{k}=\lambda,$ $k\in\mathbb{N}.$
Then, from (\ref{symR}) we have that for every $k$
\begin{equation}\label{zk}
\left\langle R\left(z_{k}, h\right)u, v\right\rangle
=\left\langle u, R\left(z_{k}, h\right)v\right\rangle.
\end{equation}

Letting $k\to\infty$ in (\ref{zk}) and using the fact that
\[R\left(\lambda+i0, h\right)=\lim_{\epsilon\downarrow 0} R(\lambda+i\epsilon, 
h) \text{ in } \mathcal{B}\left(L^{2}_{\alpha}(\mathbb{R}^{n}), 
L^{2}_{-\alpha}(\mathbb{R}^{n})\right),\,\alpha>\frac{1}{2},\] 
we obtain
\begin{equation*}
\left\langle R\left(\lambda+i0, h\right)u, v\right\rangle
=\left\langle u, R\left(\lambda+i0, h\right)v\right\rangle.
\end{equation*}
Since $C_{c}^{\infty}(\mathbb{R}^{n})\otimes C^{\infty}_{c}(\mathbb{R}^{n})$ is dense in 
$C^{\infty}_{c}(\mathbb{R}^{2n}),$ this completes the proof of the lemma.
\end{proof}

\section{Scattering Geometry}\label{sgeom}
\begin{figure}[t]
\begin{center}
\input{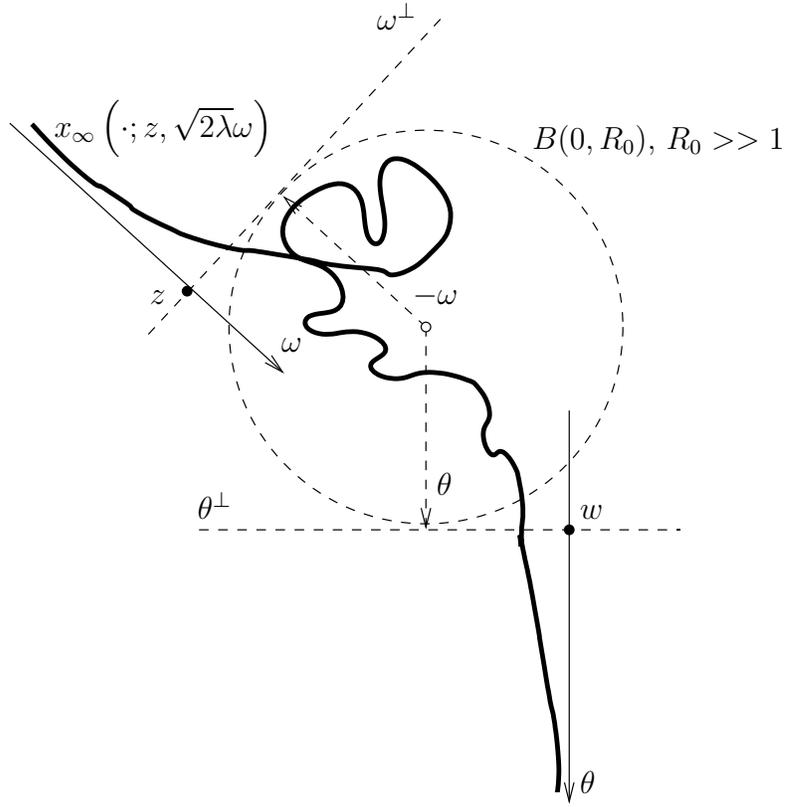}
\end{center}
\caption{The scattering relation consists of the points
$\left(\omega, z; \theta, -w\right)$ related as in this figure.}
\label{fig:sr}
\end{figure}

In this section we describe the scattering relation and prove that it can be 
parameterized by the modified actions \eqref{modaction} when the non-degeneracy assumption 
holds.
The scattering relation is a Lagrangian submanifold of $T^{*}\mathbb{S}^{n-1}\times 
T^{*}\mathbb{S}^{n-1},$ which relates the incoming and the outgoing data in the way 
suggested by Figure \ref{fig:sr}.
 
To make this precise, we first give the following 
\begin{Def}\label{NTlambda}
The trajectory $\gamma\left(\cdot; x_{0}, \xi_{0}\right)$ is non-trapped if for every $r>0$, 
there exists $T>0$ such that $\|x_{0}\|<r$ implies that
$\|x\left(s; x_{0}, \xi_{0}\right)\|>r$ for $|s|>T.$
The energy $\lambda>0$ is non-trapping if for every $r>0$ there exists $T>0$ such 
that for every $(x_{0}, \xi_{0})\in T^{*}\mathbb{R}^{n}$ with 
$\frac{1}{2}\|\xi_{0}\|^{2}+V(x_{0})=\lambda$ and $\|x_{0}\|<r$ we have $\|x\left(s; 
x_{0}, \xi_{0}\right)\|>r$ for $|s|>T.$

We also introduce the notation $T(r)$ for the infimum over $s$ 
with this property.
\end{Def}

Let, now, $\lambda>0$ be such that the operator $P(h)-\lambda$ is of principal type.
Then $\Sigma_{\lambda}=p^{-1}(\lambda)$ is a smooth $2n-1$-dimensional submanifold of $T^{*}\mathbb{R}^{n}.$

Let, further, $\left(\omega_0, z_0\right)\in T^{*}\mathbb{S}^{n-1}$ be 
such that 
$\gamma_{\infty}\left(\cdot; z_0, \sqrt{2\lambda}\omega_0\right)$ is a 
non-trapped 
trajectory with $\xi_{\infty}\left(z_0, 
\sqrt{2\lambda}\omega_0\right)\ne\omega_0.$
Then there exists $U\subset T^{*}\mathbb{S}^{n-1},$ open, $\left(\omega_0, 
z_0\right)\in U,$ such 
that for every $(\omega, z)\in U$ the trajectory $\gamma_{\infty}\left(\cdot; z, 
\sqrt{2\lambda}\omega\right)$ is non-trapped and $\xi_{\infty}\left(z, 
\sqrt{2\lambda}\omega\right)\ne\omega.$
By decreasing $U,$ if necessary, we therefore have that
\begin{equation}\label{defsr}
\begin{aligned}
SR_{\bar{U}}(\lambda)=\left\{\left(\omega, z; \xi_{\infty}\left(z, 
\sqrt{2\lambda}\omega\right), x_{\infty}\left(z, \sqrt{2\lambda}\omega\right)\right): (\omega, z)\in \bar{U}\right\}' 
\end{aligned}
\end{equation}
is a closed Lagrangian submanifold of $\left(T^{*}\mathbb{S}^{n-1}\times 
T^{*}\mathbb{S}^{n-1}, \pi_{1}^{*}\sigma+\pi_{2}^{*}\sigma\right),$ which 
we call a scattering relation at energy $\lambda$ (see Figure \ref{fig:sr}).

If $\lambda$ is a non-trapping energy level, we define the scattering relation at energy
 $\lambda>0$ as 
\begin{equation*}
SR(\lambda)=\Big\{\left(\omega, z; \xi_{\infty}\left(z, \sqrt{2\lambda}\omega\right), x_{\infty}\left(z, \sqrt{2\lambda}\omega\right)\right): \;(\omega, z)\in T^{*}\mathbb{S}^{n-1}, \omega\ne\xi_{\infty}\left(z, \sqrt{2\lambda}\omega\right)\Big\}'
\end{equation*}

We now show how, under the assumption that a certain outgoing direction is 
regular for a given incoming direction, we can find a non-degenerate phase function which 
parameterizes the scattering relation.
We begin with the following
\begin{Lem}\label{regL}
Let $\theta_0\in\mathbb{S}^{n-1}$ be regular for 
$\omega_0\in\mathbb{S}^{n-1}.$

Then there exist $O_j\subset\mathbb{S}^{n-1},$ $j=1, 2,$ open, 
$\omega_0\in O_1,$ $\theta_0\in O_2,$ and $L\in\mathbb{N}$ such 
that for every $(\omega, \theta)\in O_1\times O_2$ the number of $(\omega, 
\theta)$ trajectories is at least $L.$
\end{Lem} 

\begin{proof}
By \cite[Remark 1.1]{M} and the discussion following it, we have 
that 
there exists $L\in\mathbb{N}$ such that the number of $(\omega_0, 
\theta_0)$ trajectories is $L.$
Let $\left(z_{l}\right)_{l=1}^{L}\equiv\left(\xi_{\infty}^{-1}\left(\cdot,
\sqrt{2\lambda}\omega_0\right)\right)(\theta_0),$
By the Implicit Function Theorem, since $\theta_0$ is regular for $\omega_0,$ we have
that there exist open sets $O_1, O_2\subset\mathbb{S}^{n-1}$ with
$\theta_{0}\in O_1$ and $\omega_{0}\in O_2$ and functions $z_{l}\in
C^{\infty}(O_1\times O_2; \mathbb{R}^{n-1}),$ $l=1,
\dots, L,$ such that
$z_{l}(\omega_0, \theta_0; \lambda)=z_l$ and $\xi_{\infty}\left(z_{l}(\omega, \theta; \lambda), \sqrt{2\lambda}\omega\right)=\theta,$ $(\omega, \theta)\in O_1\times O_2,$ which completes the 
proof.
\end{proof}

Let, now, $w_{l}(\omega, \theta; \lambda)=x_{\infty}\left(z_l\left(\omega, \theta; 
\lambda\right),
\sqrt{2\lambda}\omega\right).$
As in \cite[Lemma 4]{AI}, we have the following
\begin{Lem}\label{nondeg2}
Let $\theta_0\in\mathbb{S}^{n-1}$ be regular for
$\omega_0\in\mathbb{S}^{n-1}.$

Then there exist $O_j\subset\mathbb{S}^{n-1},$ $j=1, 2,$ open,
$\omega_0\in O_1,$ $\theta_0\in O_2,$ such that the map 
\[\theta^{\perp}\ni w\mapsto\xi_{\infty}\left(w, 
-\sqrt{2\lambda}\theta\right)\in\mathbb{S}^{n-1}\] is 
non-degenerate at $w_{l}\left(\omega, \theta\right),$ $\left(\omega, \theta\right)\in O_1\times O_2,$ $l=1, \dots, L.$
\end{Lem}

We now choose $O_{1}$ and $O_{2}$ in such a way that the conclusions of Lemmas \ref{regL} 
and \ref{nondeg2} hold in some open neighborhoods of $\bar{O}_{1}$ and $\bar{O}_{2}$ and $\bar{O}_1\cap\bar{O}_2=\emptyset.$
We set 
\begin{equation}\label{lagrl}
SR_l(\lambda)=\left\{\left(\omega, \theta, z_l(\omega, \theta; \lambda), -w_l(\omega,
\theta; \lambda)
\right): (\omega, \theta)\in \bar{O}_1 \times \bar{O}_2 \right\}.
\end{equation}

The same proof as in \cite[Lemma 3.2]{RT} now shows that there exist $\bar{R}>>0,$
$T_0>T\left(\bar{R}\right),$ and open sets $U^{l}_{\omega, \theta}\subset\omega^{\perp},$ 
$z_{l}\left(\omega, \theta; \lambda\right)\in U^{l}_{\theta, \omega},$ 
$l=1, \dots, L,$ $(\omega, \theta)\in \bar{O}_1\times \bar{O}_2,$ such that
\begin{equation}\label{candefaction}
\det\left(\frac{\partial x\left(t; \cdot, 
\nabla_{x}\Phi_{-}\left(\cdot, \sqrt{2\lambda}\omega\right)(y)\right)}{\partial
y}\left(y\right)\right)\ne 0
\end{equation}
for $y\in\left\{x_{\infty}\left(s; z, 
\sqrt{2\lambda}\omega\right)\cap B\left(0, \bar{R}, \bar{R}+1\right): z\in 
U^{l}_{\omega, \theta}, s<0\right\},$ $t>T_0.$

Let, now, $t_0>T_0$ be fixed.
From \eqref{candefaction} it follows that for $(\theta, \omega)\in \bar{O}_1\times\bar{O}_2$ we define the 
(modified)
action along the segment of the $\left(\omega, 
\theta\right)$-trajectory $\gamma_l(\omega, \theta, \lambda)=\left(x_l(\omega, \theta, \lambda), \xi_l(\omega, \theta, \lambda)\right)=\gamma_{\infty}\left(\cdot; z_l(\omega, \theta; \lambda), 
\sqrt{2\lambda}\omega\right),$
between the points 
\[y_{l}\left(s; \omega, \theta, \lambda\right)=x_{\infty}\left(s; z_l\left(\omega, \theta; \lambda\right), 
\sqrt{2\lambda}\omega\right)\cap B\left(0, \bar{R}, \bar{R}+1\right)\] for 
some $s<0$ and $x_{l}(t_0; s, \omega, \theta, \lambda)
=x\left(t_0; y_l\left(s; \omega, \theta, \lambda\right), \nabla_{x}\Phi_{-}\left(y_l\left(s; 
\omega, \theta, \lambda\right), \sqrt{2\lambda}\omega\right)\right)$ and we set
\begin{equation}\label{maction}
S_{l}\left(\omega, \theta\right)=\Phi_{-}\left( y_{l}\left(s; 
\omega, \theta, \lambda\right),
\sqrt{2\lambda}\omega\right)+\int_{0}^{t_0} L\left(x, \dot{x}\right)dt-\Phi_{+}\left( x_l\left(t_0; s, \omega, \theta, 
\lambda\right),
\sqrt{2\lambda}\theta\right)+\lambda t_0,
\end{equation}
where $L(x, \dot{x})=\frac{1}{2}\left\|\dot{x}\right\|^{2}_{g}-V(x)$ is
the Lagrangian, and the integral is taken over the segment of the bicharacteristic curve $x_l(\omega, \theta, \lambda)$
connecting $y_{l}\left(s; \omega, \theta, \lambda\right)$ and 
$x_{l}\left(t_0; s, \omega, \theta, \lambda\right).$

From the representations \cite[(4.5)]{RT}
\begin{equation}\label{tail-}
\Phi_{-}(x, \sqrt{2\lambda}\omega)=2\tau\lambda+\int_{-\infty}^{\tau}\left(\frac{1}{2}\left|p_{\infty}
\left(t; z, \sqrt{2\lambda}\omega\right)\right|^{2}-V\left(q_{\infty}\left(t; z,
\sqrt{2\lambda}\omega\right)\right)-\lambda\right)dt
\end{equation}
for $x=q_{\infty}\left(\tau; z, \sqrt{2\lambda}\omega\right)\in B\left(0, \bar{R}, 
\bar{R}+1\right)$ and \cite[(4.4)]{RT}
\begin{equation}\label{tail+}
\begin{aligned}
\Phi_{+}(x, \xi)=2\lambda\tau &+\left\langle x_{\infty}\left(z, \sqrt{2\lambda}\omega\right), 
\xi\right\rangle\\
& -\int_{\tau}^{\infty}\left(\frac{1}{2}\left|p_{\infty}\left(t; z, 
\sqrt{2\lambda}\omega\right)
\right|^{2}-V\left(q_{\infty}\left(t; z, 
\sqrt{2\lambda}\omega\right)\right)-\lambda\right)dt
\end{aligned}
\end{equation}
for $(x, \xi)\in\Gamma_{+}\left(R_0, d_0, -\sigma_0\right)$ with
$x=q_{\infty}\left(\tau; z, \sqrt{2\lambda}\omega\right),$
$\xi=\lim_{t\to\infty}p_{\infty}\left(t; z, \sqrt{2\lambda}\omega\right),$ we see that $S_l(\omega, \theta)$ is independent of the choice of $s$ with the specified properties.

We now have the following 
\begin{Lem}\label{actionsr}
Let $\omega_{0}\in\mathbb{S}^{n-1}$ be regular for 
$\theta_{0}\in\mathbb{S}^{n-1}.$

Then $SR_{l}(\lambda)=\Lambda_{S_{l}},$ where $\Lambda_{S_{l}}=\left\{\left(\omega, \theta, 
d_{\omega}S_{l}, d_{\theta}S_{l}\right): (\omega, \theta)\in \bar{O}_1\times 
\bar{O}_2\right\},$ $l=1, \dots, L.$ 
\end{Lem}

\begin{proof}
We consider
\begin{equation}\label{stheta}
\begin{aligned}
d_{\theta}S_{l}(\omega, \theta) & =d_{\theta}\left(\Phi_{-}\left(
y_l\left(s; \omega, \theta, \lambda\right),
\sqrt{2\lambda}\omega\right)+\int_{0}^{t_{0}} L\left(x,
\dot{x}\right)dt\right)-d_{\omega}\Phi_{+}\left(x_l\left(t_0; s, \omega, \cdot, \lambda\right), \sqrt{2\lambda}\cdot\right)(\theta)\\
& =\left\langle\xi\left(t_0; y_l\left(s; \omega, \theta, \lambda\right), 
\nabla_{x}\Phi_{-}\left(y_{l}(s, \omega, \theta, \lambda), \sqrt{2\lambda}\omega\right)\right), d_{\theta}x_l(t_0; s, \omega, \cdot,
\lambda)(\theta)\right\rangle\\
 &\quad\quad -\left\langle\nabla_{x}\Phi_{+}\left(x_{l}(t_0; s, \omega, \theta, \lambda), \sqrt{2\lambda}\theta\right), d_{\theta}x_l(t_0; s, \omega, \cdot,
\lambda)(\theta)\right\rangle\\
& \quad\quad -d_{\theta}\left\langle \nabla_{\xi}\Phi_{+}\left(x_l(t_0; s, \omega, \theta, \lambda), \sqrt{2\lambda}\theta\right), \sqrt{2\lambda}\cdot\right\rangle(\theta)\\
& = -d_{\theta}\left\langle \nabla_{\xi}\Phi_{+}\left(x_l(t_0; s, \omega, \theta, \lambda), \sqrt{2\lambda}\theta\right), \sqrt{2\lambda}\cdot\right\rangle(\theta),
\end{aligned}
\end{equation}
where \eqref{candefaction} has allowed us to use \cite[Theorem 46.C]{A}
to obtain the second equality.
Lastly, we recall from \cite[Lemma 4.1]{RT} that 
\begin{equation}\label{xiphi+}
\lim_{t\to\infty}\left|x_{\infty}\left(t; z_l(\omega, \theta, 
\lambda),\sqrt{2\lambda}\omega\right)-\sqrt{2\lambda}\theta 
t-\nabla_{\xi}\Phi_{+}\left(x_l(t_0; s, \omega, \theta, \lambda), 
\sqrt{2\lambda}\theta\right)\right|=0.
\end{equation}

To compute $d_{\omega}S_l$ we first reparameterize the phase
trajectories in the reverse direction, which is equivalent to
considering the reverse of the initial and final directions.
Using \eqref{tail-} and \eqref{tail+} we further re-write 
$S_l\left(\omega, \theta\right)$ in the
following way
\begin{equation*}
S_l\left(\omega, \theta\right)=-\Phi_{+}\left(x_l\left(s; \omega, \theta,
\lambda\right), \sqrt{2\lambda}\theta\right) +\int_{0}^{t_{0}} L\left(x_l,
\dot{x}_l\right)dt+\Phi_{-}\left( y_l\left(t_0; s, \omega, \theta, \lambda\right),
\sqrt{2\lambda}\omega\right)+\lambda t_0,
\end{equation*}
where $x_l\left(s; \omega, \theta, \lambda\right)=x_{\infty}\left(s; z_{l}\left(\omega, \theta; \lambda\right), \sqrt{2\lambda}\omega\right)\cap B\left(0, \bar{R}, \bar{R}+1\right)$ for some $s>0,$
\[y_l\left(t_0; s, \omega, \theta, \lambda\right)=x\left(t_0;
x_l\left(s; \omega, \theta, \lambda\right),-\nabla_{x}\Phi_{+}\left(x_l\left(s; \omega, \theta, \lambda\right), \sqrt{2\lambda}\theta\right)\right),\]
and the integral is taken over the segment of the bicharacteristic curve $x_l(\omega, 
\theta, \lambda)$ connecting
$x_l\left(s; \omega, \theta, \lambda\right)$ and $y_l\left(t_0; s, \omega, \theta, \lambda\right).$
We observe that this bicharacteristic curve is uniquely defined by Lemma \ref{nondeg2} and
\eqref{candefaction}.

Lemma \ref{nondeg2} and \eqref{candefaction} further allow us to proceed as in 
(\ref{somega})
and we obtain
\begin{equation}\label{somega}
\begin{aligned}
d_{\omega}S_l\left(\omega,
\theta\right) & = d_{\omega}\left(-\Phi_{+}\left( x_l\left(s; \omega, \theta,
\lambda\right), \sqrt{2\lambda}\theta\right) +\int_{0}^{t_{0}} L\left(x,
\dot{x}\right)dt\right)\\
 &\quad +d_{\omega}\Phi_{-}\left( y_l\left(t_{0}; s, \cdot,
\theta, \lambda\right), \sqrt{2\lambda}\cdot\right)(\omega)\\
 & =d_{\omega}\left\langle \nabla_{\xi}\Phi_{-}\left(y_l(t_{0}; s, \omega, \theta, \lambda), \sqrt{2\lambda}\omega\right), \sqrt{2\lambda}\cdot\right\rangle(\omega).
\end{aligned}
\end{equation}
As above, we have that
\begin{equation}\label{xiphi-}
\lim_{t\to -\infty}\left|x_{\infty}\left(t; z_l(\omega, \theta,
\lambda),\sqrt{2\lambda}\omega\right)-\sqrt{2\lambda}\theta
t-\nabla_{\xi}\Phi_{-}\left(y_l(t_0; s, \omega, \theta, \lambda),
\sqrt{2\lambda}\omega\right)\right|=0.
\end{equation}

From (\ref{stheta}), \eqref{xiphi+}, \eqref{somega}, and \eqref{xiphi-} we 
therefore have
that $S_l$ is
a non-degenerate phase function such that $SR_l(\lambda)=\Lambda_{S_{l}}.$
\end{proof}

We remark that \eqref{tail-} and \eqref{tail+} allow us to rewrite $S_l(\omega, \theta)$ in the following way
\begin{equation}\label{laction}
\begin{aligned}
S_l(\omega, \theta)=\int_{-\infty}^{\infty}\left(\frac{1}{2}\left|p_{\infty}\left(t; z_l(\omega, \theta), \sqrt{2\lambda}\omega\right)\right|^{2}-V\left(q_{\infty}\left(t; z_l(\omega, \theta), \sqrt{2\lambda}\omega\right)\right)-\lambda\right)dt\\
-\left\langle x_{\infty}\left(z_l(\omega, \theta), \sqrt{2\lambda}\omega\right),
\sqrt{2\lambda}\theta\right\rangle,
\end{aligned}
\end{equation}
which is the same as the modified actions given by \eqref{modaction}.

\section{Proof of Main Theorem}\label{pmain}

We now turn to the proof of the Main Theorem.

\begin{proof}
Since $S(\lambda, h)$ is a unitary operator on $L^{2}(\mathbb{S}^{n-1}),$ 
we have, by Lemma \ref{tempop}, that $K_{S(\lambda, 
h)}\in\mathcal{D}'_{h}(\mathbb{S}^{n-1}\times\mathbb{S}^{n-1})$ and 
therefore $K_{T(\lambda, 
h)}\in\mathcal{D}'_{h}(\mathbb{S}^{n-1}\times\mathbb{S}^{n-1}).$

Since we are working away from the diagonal in $\mathbb{S}^{n-1}\times\mathbb{S}^{n-1}$ 
we can use integration by parts, as in \cite{RT} and \cite{M}, and obtain
\begin{equation*}
K_{T_{\pm 1}}=\mathcal{O}_{L^{2}(\mathbb{S}^{n-1}\times\mathbb{S}^{n-1}\backslash\diag 
(\mathbb{S}^{n-1}\times\mathbb{S}^{n-1}))}(h^{\infty}).
\end{equation*}
Therefore, by \eqref{kernelestim}, we obtain
\begin{equation}\label{Tpm1}
WF_{h}^{f}\left(K_{T_{\pm 1}}\right)=\emptyset.
\end{equation}

We now observe that the proof of \cite[Lemma 2.1]{RT} depends only on the estimate 
\eqref{resest} and the support properties of the symbols $a_{\pm}$ and $b_{\pm}$ and 
therefore its assertion holds here as well and we have the following 
estimates for $\gamma>\frac{n}{2}$ close to $\frac{n}{2}$
\begin{equation}\label{cutest}
\begin{aligned}
& \left\|K_{+a}^{*}(h)R(\lambda+i0, h)K_{+b}(h)\right\|_{-\gamma, 
\gamma}=\mathcal{O}(h^{\infty})\\
& \left\|K_{+a}^{*}(h)R(\lambda+i0, h)(1-\chi_b)K_{-b}(h)\right\|_{-\gamma,
\gamma}=\mathcal{O}(h^{\infty})\\
& \left\|((1-\chi_a)K_{+a})^{*}(h)R(\lambda+i0, h)K_{-b}(h)\right\|_{-\gamma,
\gamma}=\mathcal{O}(h^{\infty}),
\end{aligned}
\end{equation}
where $\chi_a\in C_{c}^{\infty}\left(B\left(0, 20R_0+1\right)\right), \chi_a(x)=1, 
|x|<20R_0$ and $\chi_b\in C_{c}^{\infty}\left(B\left(0, 10R_0+1\right)\right), \chi_b(y)=1, 
|y|<10R_0.$

From \eqref{Tpm1}, \eqref{cutest}, and \eqref{kernelestim} we then 
conclude, as in \cite[Corollary]{RT}, that
\begin{equation}\label{G0}
WF_{h}^{f}\left(\chi\left(K_{A(\lambda, h)}-c_1(n, \lambda, h)K_{G_{0}}\right)\right)=\emptyset,
\end{equation}
for every $\chi\in 
C^{\infty}\left(\mathbb{S}^{n-1}\times\mathbb{S}^{n-1}\backslash
\diag(\mathbb{S}^{n-1}\times\mathbb{S}^{n-1})\right),$ where
\begin{equation*}
G_0(\theta, \omega; \lambda, h)
=\left\langle e^{-\frac{i}{h}\Phi_{+}\left(\cdot, \sqrt{2\lambda}\theta\right)}g_{+a}(\cdot, 
\theta; h)\otimes e^{\frac{i}{h}\Phi_{-}\left(\cdot\cdot, 
\sqrt{2\lambda}\omega\right)}g_{-b}(\cdot\cdot, \omega; h), K_{R(\lambda+i0, 
h)}\right\rangle,
\end{equation*}
\begin{equation*}
g_{+a}(x, \theta; h)=e^{-\frac{i}{h}\Phi_{+}\left(x, \sqrt{2\lambda}\theta\right)}
[\chi_{a}, P_0(h)]a_{+}\left(x, 
\sqrt{2\lambda}\theta; h\right)
e^{\frac{i}{h}\Phi_{+}\left(x, \sqrt{2\lambda}\theta\right)},
\end{equation*}
\begin{equation*}
g_{-b}(y, \omega; h)=e^{-\frac{i}{h}\Phi_{-}\left(y, \sqrt{2\lambda}\omega\right)}[\chi_{b}, P_0(h)]b_{-}\left(y, \sqrt{2\lambda}\omega; h\right)
e^{\frac{i}{h}\Phi_{-}\left(y, \sqrt{2\lambda}\omega\right)},
\end{equation*}
and 
\begin{equation*}
c_1(n, \lambda, h)=2\pi(2\lambda)^{\frac{n-3}{4}}{(2\pi 
h)^{-\frac{n+1}{2}}e^{-\frac{i(n-3)\pi}{4}}}.
\end{equation*}

Let, now, $\bar{p}\in SR_{U}(\lambda)$ be such that 
$\tilde{\pi}_{1}\left(\bar{p}\right)=(\omega, z),$ where 
$\tilde{\pi}_{1}: T^{*}\mathbb{S}^{n-1}\times 
T^{*}\mathbb{S}^{n-1}\to T^{*}\mathbb{S}^{n-1}$ is the canonical projection onto the first 
factor. 
Let $A_{j}\in\Psi_{h}^{0}(1, 
\mathbb{S}^{n-1}\times\mathbb{S}^{n-1}\backslash\diag(\mathbb{S}^{n-1}\times\mathbb{S}^{n-1})),$ 
$j=0, \dots, N,$ have compactly supported symbols near $\bar{p}$ and satisfy 
$\sigma_{0}(A_j)|_{SR_{\bar{U}}(\lambda)}=0,$ $j<N.$
We also set $\varphi_{+}(x, \theta)=\Phi_{+}\left(x, \sqrt{2\lambda}\theta\right),$ 
$\left(x, 
\sqrt{2\lambda}\theta\right)\in\Gamma_{+}(R_0, d_0, \sigma_0),$ and 
$\varphi_{-}(y, \omega)=\Phi_{-}\left(y, \sqrt{2\lambda}\omega\right),$ $(y, 
\omega)\in\Gamma_{-}(R-0, d_0, -\sigma_0).$
First, we shall prove that the generalization of Egorov's 
Theorem to manifolds of unequal dimensions \cite[Lemma 7]{Afio} can be 
applied to the 
semi-classical Fourier integral operator $F$ given by the Schwartz kernel
\begin{equation*}
K_{F}=e^{-\frac{i}{h}\varphi_{+}}g_{+a}\otimes e^{\frac{i}{h}\varphi_{-}}g_{-b}.
\end{equation*}

For that, let 
\begin{equation*}
\begin{aligned}
\Lambda_{F}=\Big\{\Big( & x, y, -\nabla_{x}\varphi_{+}\left(x, \theta\right), 
\nabla_{y}\varphi_{-}\left(y, \omega\right);\, \theta, \omega, 
-\nabla_{\theta}\varphi_{+}\left(x, \theta\right), \nabla_{\omega}\varphi_{+}\left(y, \omega\right):\\
 & \left(x, \sqrt{2\lambda}\theta\right)\in\Gamma_{+}(R_0, d_0, \sigma_0)\cap 
\left(T^{*}(\supp\nabla\chi_{a})\times T^{*}\mathbb{S}^{n-1}_{2\lambda}\right),\\
& \left(y, \sqrt{2\lambda}\omega\right)\in\Gamma_{-}(R_0, d_0, 
-\sigma_0)\cap\left(T^{*}(\supp\nabla\chi_{b})\times T^{*}\mathbb{S}^{n-1}_{2\lambda}\right)\Big)\Big\}.
\end{aligned}
\end{equation*}
For every 
$(x, \xi)\in\Gamma_{\pm}(R_0, d_0, \pm\sigma_0)$ there exist unique phase trajectories
$(q_{\pm}(\cdot; x, \xi), p_{\pm}(\cdot; x, \xi))$ such that $q_{\pm}(0; x, \xi)=x$ and 
$\lim_{t\to\pm\infty}p_{\pm}(t; x, \xi)=\xi,$ respectively (see \cite[Subsection 4.1]{RT} as 
well as the discussion following \cite[Definition 1.10]{I}).
Furthermore, by the construction of $\Phi_{\pm},$
\begin{equation*}
\nabla_{x}\Phi_{\pm}(q_{\pm}(t; x, \xi), \xi)=p_{\pm}(t; x, \xi).
\end{equation*}
By \cite[Lemma 4.1]{RT}, we also have that 
\begin{equation*}
\lim_{t\to\pm\infty}\left|q_{\pm}(t; x, \xi)-\xi t-\nabla_\xi\Phi_\pm(x, \xi)\right|=0.
\end{equation*}
These considerations imply that 
\begin{equation}\label{immerse}
\pi_{1}|_{\Lambda_{F}} \text{ is an immersion,}
\end{equation}
where 
\begin{equation*}
\pi_{1}:T^{*}\mathbb{R}^{n}\times T^{*}\mathbb{R}^{n}\times T^{*}\mathbb{S}^{n-1}\times T^{*}\mathbb{S}^{n-1}\to 
T^{*}\mathbb{R}^{n}\times T^{*}\mathbb{R}^{n}
\end{equation*}
is the canonical projection.
With \eqref{immerse} the hypotheses of the generalization of 
Egorov's 
Theorem to manifolds of unequal dimensions \cite[Lemma 7]{Afio} are 
satisfied and applying \cite[Lemma 7]{Afio} we obtain that there exist 
$B_{j}\in\Psi_{h}^{0}(1, 
\mathbb{R}^{n}\times\mathbb{R}^{n}),$ $j=0, \dots, N,$ satisfying the 
following conditions
\begin{enumerate}
\item $\sigma(B_j),$ $j=0, \dots, N,$ have compact support near 
a point $\bar{q}\in T^{*}\mathbb{R}^{n}\times T^{*}\mathbb{R}^{n}$ such that 
$\hat{\pi}_{1}\left(\bar{q}\right)\in\gamma_{\infty}\left(\cdot; z, \sqrt{2\lambda}\theta\right),$ where 
$\hat{\pi}_{1}: T^{*}\mathbb{R}^{n}\times T^{*}\mathbb{R}^{n}\to T^{*}\mathbb{R}^{n}$ 
is the canonical projection onto the first factor.
\item  $\sigma_{0}(B_j)|_{\Lambda_{R}(\lambda)}=0, j<N,$ where 
$\Lambda_R(\lambda)=\cup_{t>0}\left(\graph\exp(tH_{p})|_{\Sigma_{\lambda}}\right)'.$
\item Near $\left(\bar{p}, \bar{q}\right),$ 
\begin{equation}\label{interwine}
\left(\prod_{j=0}^{N}A_{j}\right)\left(e^{-\frac{i}{h}\phi_{+}}g_{+a}\otimes 
e^{\frac{i}{h}\phi_{-}}g_{-b}\right)\equiv 
\left(e^{-\frac{i}{h}\phi_{+}}g_{+a}\otimes e^{\frac{i}{h}\phi_{-}}g_{-b}\right)
\left(\prod_{j=0}^{N}B_j\right).
\end{equation} 
\end{enumerate}

Assumption \eqref{resest} and Lemma \ref{tempop}, now, imply that 
$K_{R(\lambda+i0, h)}\in\mathcal{D}_{h}'(\mathbb{R}^{2n}).$ 
From \eqref{interwine} we therefore obtain
\begin{equation}\label{interw2}
\left(\prod_{j=0}^{N}A_{j}\right)K_{A(\lambda, h)}\equiv 
c_{1}\left(n, \lambda, h\right)
\left(e^{-\frac{i}{h}\phi_{+}}g_{+a}\otimes e^{\frac{i}{h}\phi_{-}}g_{-b}\right)
\left(\prod_{j=0}^{N}B_j\right)\left(\chi_{2}\otimes\chi_{1}\right)K_{R(\lambda+i0, h)},
\end{equation}
near $\left(\bar{p}, \bar{q}\right),$ where $\chi_j\in 
C_{c}^{\infty}(\mathbb{R}^{n}; \mathbb{R}),$ $j=1, 2,$ are such that
$\chi_2=1$ on $\supp g_{+a},$ $\chi_1=1$ on $\supp g_{-b},$ and
$\supp\chi_1\cap\supp\chi_2=\emptyset.$

Estimate \eqref{resest}, Lemma \ref{rsym}, and the same proof as in 
\cite[Theorem 1]{AI} further give that there exists an open set 
$V\subset\Lambda_{R}(\lambda),$ $\bar{q}\in V,$ such that
$\left(\chi_{2}\otimes\chi_{1}\right)K_{R(\lambda+i0, h)}\in 
I_{h}^{1}\left(\mathbb{R}^{2n}, \Lambda_{R}(\lambda)\cap\bar{V}\right).$ 
(We recall here that the fact that $\chi_1$ and $\chi_2$ have disjoint support is crucial in 
the proof of \cite[Theorem 1]{AI}.)
Therefore
\begin{equation}\label{estR}
\left(\prod_{j=0}^{N}B_j\right)\left(\chi_{2}\otimes\chi_{1}\right)K_{R(\lambda+i0, 
h)}
=\mathcal{O}_{L^{2}(\mathbb{R}^{2n})}\left(h^{N-1-\frac{n}{2}}\right), h\to 0.
\end{equation}
Since $g_{+b}, g_{-a}\in S^{-1}_{2n-1}(1)\cap 
C_{c}^{\infty}(\mathbb{R}^{n}\times\mathbb{S}^{n-1}),$ we easily find that
\begin{equation}\label{estF}
\left\|\left(e^{-\frac{i}{h}\phi_{+}}g_{+a}\otimes e^{\frac{i}{h}\phi_{-}}g_{-b}\right) 
\right\|_{\mathcal{B}(L^{2}(\mathbb{R}^{n}), L^{2}(\mathbb{S}^{n-1}))}=\mathcal{O}(h). 
\end{equation}
Estimates \eqref{estR} and \eqref{estF} together with \eqref{G0} and 
\eqref{interw2} now imply that
\begin{equation*}
\left(\prod_{j=0}^{N}A_{j}\right)K_{A(\lambda, h)}=
\mathcal{O}_{L^{2}(\mathbb{S}^{n-1}\times\mathbb{S}^{n-1})}\left(h^{N-n-\frac{3}{2}}\right),
\end{equation*}
and therefore
\begin{equation*}
A(\lambda, h)\in \mathcal{I}_{h}^{\frac{n}{2}+2}\left({\mathbb
S}^{n-1} \times{\mathbb S}^{n-1}\backslash\diag(\mathbb{S}^{n-1}\times\mathbb{S}^{n-1}),
SR_{\bar{U}}(\lambda)\right).\qedhere
\end{equation*}
\end{proof}

\section{Applications}
In this section we discuss two applications of our Main Theorem to trapping and non-trapping 
energies, respectively. 
\subsection{Non-Trapping Energies}\label{snt}
\begin{Co}\label{ntfio}
Let $\lambda>0$ be a non-trapping energy level for $P$ and such that $P(h)-\lambda$ is of 
principal type.

Then $A(\lambda, h)\in \mathcal{I}_{h}^{\frac{n}{2}+2}\left({\mathbb
S}^{n-1} \times{\mathbb S}^{n-1}\backslash\diag(\mathbb{S}^{n-1}\times\mathbb{S}^{n-1}),
SR(\lambda)\right).$
\end{Co}

\begin{proof} 
From \cite[Lemma 2.2]{RT} we have that $\left\|R(\lambda+i0, h)\right\|_{\alpha, 
-\alpha}=\mathcal{O}\left(\frac{1}{h}\right), \alpha>\frac{1}{2}.$
The result now follows from the Main Theorem.
\end{proof}

\subsection{Trapping Energies}\label{str}
\begin{Co}\label{tfio}
Let $\lambda>0$ be a trapping energy level for $P$ and such that $P(h)-\lambda$ is of
principal type.
Let also
\begin{enumerate}[(i)]
\item there exist $\theta_0\in [0, \pi),$ $R>0$ such that the potential $V$ 
extends holomorphically to the 
domain $D_{R, \theta_{0}}=\{z\in\mathbb{C}^{n}: |z|>R, |\Im z|\leq \tan\theta_0|\Re z|\}$ and $|V(x)|\leq C|x|^{-\beta}$ for all $x\in D_{R, \theta_{0}}$ and some $\beta>0,$ $C>0,$
and
\item $Res(P(h))\cap([\lambda-\epsilon, \lambda+\epsilon]+i[0, Ch^{M}])=\emptyset$ for some $\epsilon>0,$ $C>0,$ and $M>0.$  
\end{enumerate}
Lastly, let there exist $(\theta, z)\in T^{*}\mathbb{S}^{n-1}$ such that 
$\gamma_{\infty}\left(\cdot; z, \sqrt{2\lambda}\theta\right)$ is a 
non-trapped 
trajectory.

Then there exists an open set $U\subset T^{*}\mathbb{S}^{n-1}$ such that 
\begin{equation*}
A(\lambda, h)\in 
\mathcal{I}_{h}^{\frac{n}{2}+2}\left({\mathbb S}^{n-1} \times{\mathbb 
S}^{n-1}\backslash\diag(\mathbb{S}^{n-1}\times\mathbb{S}^{n-1}),
SR_{\bar{U}}(\lambda)\right).
\end{equation*}
\end{Co}

\begin{proof} 
We choose $U$ as in Definition \ref{defsr}.
From \cite[Proposition 4.1]{M}, we have that 
there exists 
$m\in\mathbb{N}$ such that 
\begin{equation*}
\left\|R(\lambda+i0, h)\right\|_{\alpha, -\alpha}
=\mathcal{O}\left(\frac{1}{h^{m}}\right), \alpha>\frac{1}{2}.
\end{equation*}
The assertion of the Corollary now follows from the Main 
Theorem.
\end{proof}

\subsection{Microlocal Representation of the Scattering Amplitude}\label{smicrol}
Here we show how under the non-degeneracy assumption the expansion \eqref{vexpansion}
follows from the results we have proved in this article and the
characterization of semi-classical
 Fourier integral distributions as oscillatory integrals, which we have
developed in \cite[Theorem 1]{AI}.
More precisely, we have the following

\begin{Th}\label{tmicrol}
Let $\omega_0\in\mathbb{S}^{n-1}$ be regular for
$\theta_0\in\mathbb{S}^{n-1}$ and
$L\in\mathbb{N}$ be the number of $(\theta_0, \omega_0)$ phase
trajectories.
Let $\lambda>0$ be such that $P-\lambda$ is of principle type and 
$\left\|R(\lambda+i0, h)\right\|_{\alpha, -\alpha}=\mathcal{O}(h^{m}),$ 
$m\in\mathbb{R},$ $\alpha>\frac{1}{2}.$

Then, if $P_l\in\Psi_{h}^{0}(1, \mathbb{S}^{n-1}\times\mathbb{S}^{n-1}),$ $l=1,
\dots, L,$ are microlocal cut-offs to the Lagrangian submanifolds
$SR_l(\lambda)$ defined by
\eqref{lagrl},
respectively, 
\begin{equation*}
P_l K_{A(\lambda, h)}=e^{\frac{i}{h}S_{l}}a_{l},\: l=1, \dots, L,
\end{equation*}
where $S_{l},$ $l=1, \dots, L,$ are as given by
(\ref{maction}) and $a_l\in S_{2n-2}^{n+\frac{3}{2}}(1),$ $l=1, \dots, L,$
have compact support.
\end{Th}

\begin{proof}
By our Main Theorem, $A(\lambda, 
h)\in\mathcal{I}_{h}^{\frac{n}{2}+2}\left(\mathbb{S}^{n-1}\times\mathbb{S}^{n-1}, 
\cup_{l=1}^{L}SR_{l}(\lambda)\right).$
From \cite[Lemma 5]{AI} it follows that $P_l K_{A(\lambda,
h)}\in I_{h}^{\frac{3}{2}}\left(\mathbb{S}^{n-1}\times\mathbb{S}^{n-1},
SR_{l}(\lambda)\right),$ $l=1, \dots, L.$
With this and Lemma \ref{actionsr} the hypotheses of \cite[Theorem 1]{AI} are
satisfied and we obtain that there exist $a_l\in S_{2n-2}^{n+\frac{3}{2}}(1),$ $l=1, \dots,
L,$ such that $P_l K_{A(\lambda, h)}=e^{\frac{i}{h}S_{l}}a_{l}$ microlocally near
$SR_{l}(\lambda),$ $l=1, \dots, L.$
\end{proof}

We remark that the conclusion of this theorem holds whenever we have polynomial bound on the 
resolvent.
We also remark that this theorem recovers the phases \eqref{modaction} in 
\eqref{vexpansion}, due to \eqref{laction}.

\appendix
\section{Elements of Semi-Classical Analysis}\label{scanal}
In this section we recall some of the elements of semi-classical analysis
which we use in this paper.
First we define two classes of symbols
\begin{equation*}
S_{2n}^{m}\left(1\right)= \left\{ a\in
C^{\infty}\left(\mathbb{R}^{2n}\times(0, h_0]\right): \forall
\alpha, \beta\in\mathbb{N}^{n}, \sup_{(x, \xi,
h)\in\mathbb{R}^{2n}\times (0,
h_{0}]}h^{m}\left|\partial^{\alpha}_{x}\partial^{\beta}_{\xi}a\left(x,
\xi;
h\right)\right|\leq
C_{\alpha, \beta}\right\}
\end{equation*}
and
\begin{equation*}
S^{m, k}\left(T^{*}\mathbb{R}^{n}\right)=\left\{a\in
C^{\infty}\left(T^{*}\mathbb{R}^{n}\times(0, h_0]\right): \forall \alpha,
\beta\in\mathbb{N}^{n}, \left|\partial^{\alpha}_{x}\partial^{\beta}_{\xi}
a\left(x, \xi;
h\right)\right|\leq
C_{\alpha,
\beta}h^{-m}\left\langle\xi\right\rangle^{k-|\beta|}\right\},
\end{equation*}
where $h_0\in(0,1]$ and $m, k\in\mathbb{R}.$
For $a\in S_{2n}\left(1\right)$ or $a\in S^{m, 
k}\left(T^{*}\mathbb{R}^{n}\right)$ we define
the
corresponding semi-classical pseudodifferential operator of class
$\Psi_{h}^{m}(1, \mathbb{R}^{n})$ or $\Psi_{h}^{m, k}(\mathbb{R}^{n}),$
respectively, by
setting
\begin{equation*}
Op_{h}\left(a\right)u\left(x\right)=\frac{1}{\left(2\pi
h\right)^{n}}\int\int e^{\frac{i\left\langle x-y,
\xi\right\rangle}{h}}a\left(x, \xi; h\right)u\left(y\right) dy d\xi, 
\;u\in
\mathcal{S}\left(\mathbb{R}^{n}\right),
\end{equation*}
and extending the definition to $\mathcal{S}'\left(\mathbb{R}^{n}\right)$ 
by
duality (see \cite{DS}).
Here we work only with symbols which admit asymptotic expansions in 
$h$ and with
pseudodifferential operators which are quantizations of such symbols.
For $A\in\Psi_{h}^{k}(1, \mathbb{R}^{n})$ or $A\in\Psi_{h}^{m, 
k}(\mathbb{R}^{n}),$ we shall
use $\sigma_{0}(A)$ and $\sigma(A)$ to denote its principal symbol and its 
complete symbol,
respectively.
A semi-classical pseudodifferential operator is said to be of principal
type if its
principal symbol $a_0$ satisfies
\begin{equation}\label{prtype}
a_0=0\implies da_0\ne 0.
\end{equation}

For $a\in S^{m, k}_{n}\left(T^{*}\mathbb{R}^{n}\right)$ we define:
\begin{equation*}
\begin{aligned}
& \esupp_{h} a\\
&\quad =\Big\{\left(x, \xi\right)\in T^{*}\mathbb{R}^{n}|\: \exists\:
\epsilon>0\;
\partial_{x}^{\alpha}\partial_{\xi}^{\beta}a\left(x',
\xi'\right)=\mathcal{O}_{C(B((x, \xi),
\epsilon))}\left(h^{\infty}\right),\; \forall\alpha,
\beta\in\mathbb{N}^{n}\Big\}^{c}\\
& \quad\cup\bigg(\bigg\{\left(x, \xi\right)\in
T^{*}\mathbb{R}^{n}\backslash\left\{0\right\}|\:
\exists\:\epsilon>0\:
\partial_{x}^{\alpha}\partial_{\xi}^{\beta}a\left(x',
\xi'\right)=\mathcal{O}\left(h^{\infty}\left\langle\xi\right\rangle^{-\infty}\right),\\
&\quad\quad \text{uniformly in } (x', \xi') \text{ such that }
\|x-x'\|+\frac{1}{\|\xi'\|}+\left\|\frac{\xi}{\|\xi\|}-\frac{\xi'}{\|\xi'\|}\right\|
<\epsilon\bigg\}/
\mathbb{R}_{+}\bigg)^{c}\\
&\quad \subset T^{*}\mathbb{R}^{n}\sqcup S^{*}\mathbb{R}^{n},
\end{aligned}
\end{equation*}
where we define
$S^{*}\mathbb{R}^{n}=\left(T^{*}\mathbb{R}^{n}\backslash\left\{0\right\}\right)/\mathbb{R}_{+}.$
For $A\in\Psi^{m, k}_{h}\left(\mathbb{R}^{n}\right),$ we then define
\begin{equation*}
WF_{h}\left(A\right)=\esupp_{h} a, A=Op_{h}\left(a\right).
\end{equation*}

We also define the class of semi-classical distributions
$\mathcal{D}_{h}'(\mathbb{R}^{n})$ with which we will work here
\begin{equation*}
\begin{aligned}
\mathcal{D}'_{h}(\mathbb{R}^{n}) = & \big\{u\in C^{\infty}_{h}\left((0,
1];
\mathcal{D}'\left(\mathbb{R}^{n}\right)\right): \forall\chi\in
C_{c}^{\infty}\left(\mathbb{R}^{n}\right) \exists\: N\in\mathbb{N}\text{
and
} C_{N}>0:\\
& \quad |\mathcal{F}_{h}\left(\chi u\right)\left(\xi, h\right)|\leq
C_{N}h^{-N}\langle\xi\rangle^{N}\big\}
\end{aligned}
\end{equation*}
where
\begin{equation*}
\mathcal{F}_{h}\left(u\right)\left(\xi,
h\right)=\int_{\mathbb{R}^{n}}e^{-\frac{i}{h}\left\langle x,
\xi\right\rangle}u\left(x, h\right)dx
\end{equation*}
with the obvious extension of this definition to
$\mathcal{E}_{h}'(\mathbb{R}^{n}).$
We work with the $L^{2}-$based semi-classical Sobolev spaces 
$H^{s}(\mathbb{R}^{n}),$ $s\in\mathbb{R},$ which consist of the 
distributions   
$u\in\mathcal{D}_{h}'(\mathbb{R}^{n})$ such that 
$\|u\|_{H^{s}(\mathbb{R}^{n})}^{2}\overset{\df}{=}\frac{1}{(2\pi 
h)^{n}}\int_{\mathbb{R}^{n}}(1+\|\xi\|^{2})^{s}\left|\mathcal{F}_{h}(u)(\xi, 
h)\right|^{2}d\xi<\infty.$ 

For $u\in\mathcal{D}_{h}'(\mathbb{R}^{n})$ we also define its finite semi-classical 
wavefront set as follows.
\begin{Def}\label{defM}
Let $u\in\mathcal{D}'_{h}\left(\mathbb{R}^{n}\right)$ and let
$\left(x_{0}, \xi_{0}\right)\in\hat{T}^{*}\left(\mathbb{R}^{n}\right).$
Then the point $\left(x_{0}, \xi_{0}\right)$ does not belong to
$WF_{h}^{f}\left(u\right)$ if there exist $\chi\in
C_{c}^{\infty}\left(\mathbb{R}^{n}\right)$ with $\chi\left(x_{0}\right)\ne
0$ and an open neighborhood
$U$ of $\xi_{0}$, such that $\forall N\in\mathbb{N},$ $\forall\xi\in U,$
$|\mathcal{F}\left(\chi u\right)\left(\xi, h\right)|\leq C_{N}h^{N}.$
\end{Def}

We say that $u=v$ {\it microlocally} (or $u\equiv v$) near an open set
$U\subset T^{*}\mathbb{R}^{n}$, if
$P(u-v)=\mathcal{O}\left(h^{\infty}\right)$ in
$C_{c}^{\infty}\left(\mathbb{R}^{n}\right)$ for
every $P\in \Psi^{0}_{h}\left(1, \mathbb{R}^{n}\right)$ such that
\begin{equation}\label{P}
WF_{h}\left(P\right)\subset \tilde{U}, \bar{U}\Subset \tilde{U}\Subset
T^{*}\mathbb{R}^{n}, \tilde{U} \text{ open}.
\end{equation}
We also say that $u$ satisfies a property $\mathcal{P}$  {\it
microlocally} near an open set $U\subset T^{*}{\mathbb{R}^{n}}$ if there
exists $v\in\mathcal{D}_{h}'\left(\mathbb{R}^{n}\right)$ such that $u=v$
microlocally near $U$ and $v$ satisfies property $\mathcal{P}$.

For open sets $U, V\subset T^{*}\mathbb{R}^{n},$ the operators $T,
T'\in\Psi^{m}_{h}\left(\mathbb{R}^{n}\right)$ are said to be {\it
microlocally
equivalent} near $V\times U$ if for any $A,
B\in\Psi_{h}^{0}\left(\mathbb{R}^{n}\right)$
such that
\begin{equation*}
WF_{h}\left(A\right)\subset\tilde{V},
WF_{h}\left(B\right)\subset\tilde{U},
\bar{V}\Subset\tilde{V}\Subset T^{*}\mathbb{R}^{n},
\bar{U}\Subset\tilde{U}\Subset T^{*}\mathbb{R}^{n}, \tilde{U}, \tilde{V}
\text{ open }
\end{equation*}
\begin{equation*}
A\left(T-T'\right)B=\mathcal{O}\left(h^{\infty}\right)\colon\mathcal{D}_{h}'
\left(\mathbb{R}^{n}\right)\rightarrow C^{\infty}\left(\mathbb{R}^{n}\right).
\end{equation*}
We also use the notation $T\equiv T'.$

We extend these notions to compact manifolds through the following definition of a
semi-classical pseudodifferential operator on a compact manifold.
Let $M$ be a smooth compact manifold and $\kappa_{j}: M_j\to X_j,$ $j=1, \dots, N,$ a set of
local charts.
A linear continuous operator $A: C^{\infty}(M)\to \mathcal{D}_{h}'(M)$ belongs to
$\Psi_{h}^{m}(1, M)$ if for all $j\in\{1, \dots, N\}$ and $u\in C^{\infty}_{c}(M_{j})$ we 
have
$Au\circ\kappa_{j}^{-1}=A_j\left(u\circ\kappa_{j}^{-1}\right)$ with
$A_j\in\Psi_{h}^{m}\left(X_j\right),$ and
$\chi_1 A\chi_2:\mathcal{D}_{h}'(M)\to h^{\infty}C^{\infty}(M)$ for $\chi_j\in
C^{\infty}(M)$ with
$\supp\chi_1\cap\supp\chi_2=\emptyset.$

Lastly, we define global semi-classical Fourier integral operators.
\begin{Def}\label{dfio}
Let $M$ be a smooth $k$-dimensional manifold and let $\Lambda\subset
T^{*}M$ be a smooth closed
Lagrangian submanifold with respect to the canonical symplectic
structure on $T^{*}M.$
Let $r\in\mathbb{R}.$
Then the space $I^{r}_{h}\left(M, \Lambda\right)$ of semi-classical
Fourier integral
distributions of order $r$ associated to $\Lambda$ is defined as the set
of all $u\in\mathcal{D}'_{h}\left(M\right)$
such
that
\begin{equation}\label{defgfio}
\left(\prod_{j=0}^{N}
A_{j}\right)\left(u\right)=\mathcal{O}_{L^{2}\left(M\right)}\left(h^{N-r-\frac{k}{4}}\right),
h\to 0,
\end{equation}
for all $N\in\mathbb{N}_{0}$ and for all $A_{j}\in \Psi_{h}^{0}\left(1,
M\right),$ $j=0, \dots, N-1,$ with
compactly
supported symbols and principal symbols vanishing on $\Lambda$, and any $
A_N \in
\Psi_h^{0} ( 1 , M ) $ with a compactly supported symbol.

A continuous linear operator
$C_{c}^{\infty}\left(M_1\right)\rightarrow\mathcal{D}_{h}'\left(M_2\right),$
where $M_1, M_2$ are smooth manifolds,
whose Schwartz kernel is an element of
$I_{h}^{r}(M_1\times M_2, \Lambda)$ for some
Lagrangian submanifold $\Lambda\subset T^{*}M_1\times T^{*}M_2$ and some $r\in\mathbb{R}$
will be called a global semi-classical Fourier integral
operator of order $r$ associated to $\Lambda.$
We denote the space of these operators by
$\mathcal{I}_{h}^{r}(M_1\times M_2, \Lambda).$
\end{Def}

{\bf Acknowledgements.}  I would like to thank Maciej Zworski for the idea to use the proof 
of Schwartz Kernel Theorem in the proof of Lemma \ref{tempop}.

\end{document}